\input amstex 
\documentstyle{amsppt} 
\magnification=1200 
 
\loadbold
 
\TagsOnRight 
\NoRunningHeads 
\NoBlackBoxes

\def\la{\lambda} 
\def\La{\Lambda} 
\def\om{\omega} 
\def\ep{\varepsilon} 
\def\de{\delta} 
 
\def\C{\Bbb C}

\def\Z{\Bbb Z}

\def\dim{\operatorname{dim}} 
\def\Tab{\operatorname{Tab}}

\def\wha{\widehat a} 
\def\whb{\widehat b} 
\def\half{\tfrac12}

\def\const{\operatorname{const}} 
\def\T{\Cal T}  
\def\sgn{\operatorname{sgn}} 

\def\bp{\boldkey p}

\topmatter 
\title Frobenius--Schur functions: summary of results \endtitle 
\author Grigori Olshanski, Amitai Regev, and Anatoly Vershik  
\endauthor 

\endtopmatter

\document

We introduce and study a family $\{FS_\mu\}$ of symmetric functions
which we call the {\it Frobenius--Schur functions.\/} These are
inhomogeneous functions indexed by partitions and such that $FS_\mu$
differs from the conventional Schur function $s_\mu$ in lower terms
only. 

Our interest in these new functions comes from the fact that they provide
an explicit expression for $\dim\nu/\mu$ (dimension of a skew Young
diagram $\nu/\mu$) in terms of the Frobenius coordinates of $\nu$,
see Theorem 2 below. 
The problem of studying $\dim\nu/\mu$ as a function of the Frobenius
coordinates of $\nu$ originated in the asymptotic theory of the
characters of the symmetric groups \cite{VK}. We discuss this point below.

Our main result is a surprisingly simple determinantal expression of
the Frobenius--Schur functions in terms of the conventional Schur
functions, see Theorem 9.

Other results concerning the functions $FS_\mu$ are: certain
generating series, the Giambelli 
formula, vanishing properties and interpolation, a combinatorial
formula (representation in terms of tableaux), and a
Sergeev--Pragacz--type  formula. 

We actually deal with a wider class of symmetric
functions $s_{\mu;a}$ of the form $s_\mu+\text{lower terms}$, which
we call {\it multiparameter Schur functions.\/} Here
$a=(a_i)_{i\in\Z}$ is an arbitrary doubly infinite sequence of
parameters. The multiparameter Schur functions turn into the
Frobenius--Schur functions when $(\ldots, a_0,a_1,a_2,a_3,\ldots)
=(\ldots,-\frac12,\frac12,\frac32,\frac52,\ldots)$ and into the
conventional Schur functions when $a_i\equiv0$. All our results about
the functions $FS_\mu$ can be stated for the functions
$s_{\mu;a}$ without extra efforts. For this reason we decided to include
this generalization (although at present we do not dispose of any
applications of the functions $s_{\mu;a}$ with general parameters). 

For proofs we refer to the detailed version of the paper, \cite{ORV}. 

We are grateful to Alain Lascoux for fruitful criticism.

The authors were supported by Russian Foundation
for Basic Research under grants 98--01--00303 and 99--01--00098 (G.~O. and
A.~V., respectively) and by ISF grant 6629/1 (A.~R.). 

\subhead Preliminaries \endsubhead 
We follow the notation of Macdonald
\cite{Ma1} in what concerns symmetric functions and Schur functions.
Our base field is $\C$. Let $\La$ denote the algebra of symmetric
functions. Recall that 
$\La$ is a graded algebra, isomorphic to the algebra of polynomials 
in the power sums $\bp_1,\bp_2,\dots\,$. A natural homogeneous basis in 
$\La$ is formed by the {\it Schur functions\/} $s_\mu$, where the
subscript $\mu$ is a partition (or, what is the same, a Young
diagram). 

The representation theory of the infinite symmetric group naturally
leads to linear superalgebra and supersymmetric functions, see
\cite{VK}, \cite{W}, \cite{KV1}, \cite{KV2}, \cite{Ol}, \cite{KOO}. 

A polynomial in two finite collections of indeterminates,
$x_1,\dots,x_m$ and $y_1,\dots,y_n$, is called {\it supersymmetric\/}
if it is separately symmetric in the $x$'s and in the $y$'s, and,
moreover, the following {\it cancellation property\/} holds: for any
fixed $i$ and $j$, specializing $x_i=-y_j=t$ we get a polynomial in
$m+n-2$ variables not depending on $t$.  

Let $\C^\infty_0$ be the space of the infinite complex vectors with finitely
many nonzero coordinates. A function $f(x;y)$ on
$\C^\infty_0\times\C^\infty_0$ is called a {\it supersymmetric
function\/} if for any $m$ and $n$, the specialization
$x_{m+1}=x_{m+2}=\dots=y_{n+1}=y_{n+2}=\dots=0$ turns $f$ into a
supersymmetric polynomial whose degree remains bounded as $m$ and $n$ go to
infinity. 

Specializing $\bp_k\mapsto \sum_i x_i^k+(-1)^{k-1}\sum_j y_j^k$ we
get an {\it isomorphism\/} between the algebra $\La$ and the algebra
of supersymmetric functions, see \cite{S}, \cite{P}. We agree to identify both
algebras by 
means of this isomorphism. Given $f\in\La$, we denote by $f(x;y)$ the
corresponding supersymmetric function on
$\C^\infty_0\times\C^\infty_0$. In particular, we write $s_\mu(x;y)$
for the ``supersymmetric Schur function''. 

If $x\in\C^m$, $y\in\C^n$, and $f$ is a supersymmetric function, then
we abbreviate $f(x;y)=f(x,0,0,\dots;y,0,0,\dots)$. 

We will often use the {\it Frobenius notation\/} \cite{Ma1, \S I.1}
for Young diagrams: 
$$ 
\nu=(p_1,\dots,p_d\mid q_1,\dots,q_d),   
$$
where $d=d(\nu)$ (the {\it depth\/} of $\nu$) is the number of 
diagonal squares in $\nu$, $p_i=\nu_i-i$ and $q_i=\nu'_i-i$ are the
{\it Frobenius coordinates\/} of $\nu$, and $\nu'$ stands for the
transposed diagram.  
As usual, $|\nu|$ denotes the number of squares in $\nu$.

\subhead The Frobenius--Schur functions: a motivation\endsubhead
Let $S(n)$ denote the symmetric group of degree $n$. For a Young
diagram $\nu$ let $\chi^\nu$ denote the irreducible character of
$S(|\nu|)$ indexed by $\nu$ and let $\dim\nu=\chi^\nu(e)$ be the
degree of $\chi^\nu$ (i.e., the dimension of the corresponding
representation), which coincides with the number 
of standard tableaux of shape $\nu$. If $\sigma$ is a partition of
$|\nu|$ then $\chi^\nu_\sigma$ denotes the value of $\chi^\nu$ on the
conjugacy class indexed by $\sigma$. Our starting point is the
following result (see \cite{VK}, \cite{W, III.6}, \cite{KO}):

\proclaim{Theorem 1} Fix a partition $\rho=(\rho_1,\dots,\rho_l)$ 
of $m$. For any $n\ge m$ we agree to complete $\rho$ to a partition of $n$
by adding 1's. Let $\nu=(p_1,\dots,p_d\mid q_1,\dots,q_d)$ denote an
arbitrary Young diagram with $n$ boxes.  

There
exists an inhomogeneous symmetric function 
$$
\bp^\#_\rho=\bp_{\rho_1}\dots\bp_{\rho_l}\, +\, \text{\rm lower degree
terms} 
$$
such that 
$$
\frac{\chi^\nu_{\rho\cup 1^{n-m}}}{\dim\nu}=
\frac{\bp^\#_\rho(p_1+\frac12,\dots,p_d+\frac12;
q_1+\frac12,\dots,q_d+\frac12)}
{n(n-1)\dots(n-m+1)}\,.
$$
\endproclaim
To avoid misunderstanding, let us explain once more the meaning of an
expresion like $\bp^\#_\rho(p_1+\frac12,\dots,p_d+\frac12;
q_1+\frac12,\dots,q_d+\frac12)$: here a symmetric function
(specifically, $\bp^\#_\rho$) is first viewed as a supersymmetric
function on $\C^\infty_0\times\C^\infty_0$ and then it is specialized
at a finite vector from the space $\C^d\times\C^d$ canonically
embedded into the larger space $\C^\infty_0\times\C^\infty_0$. 
 
It is worth noting
that a symmetric function is uniquely determined by its
superspecializations at the points of the form 
$(p_1+\frac12,\dots,p_d+\frac12;q_1+\frac12,\dots,q_d+\frac12)$, so
that $\bp^\#_\rho$ is defined uniquely.

Note also a discrepancy with the classical point of view; usually a
character $\chi^\nu$ is thought of as a function on the set of conjugacy
classes whereas here the picture is reversed: we fix a conjugacy class
and view $\chi^\nu$ as a function of the index $\nu$. 

This theorem plays an important role in the following {\it asymptotic
problem\/} considered in \cite{VK}: describe all sequences $\{\nu\}$
of Young diagrams (where $|\nu|\to\infty$) such that the corresponding
normalized irreducible characters $\chi^\nu/\dim\nu$ converge on the
infinite symmetric group $S(\infty)$ (this group is defined as
the union of the groups $S(n)$). It turns out that the limit
functions $\lim\limits_{\{\nu\}}\frac{\chi^\nu}{\dim\nu}$ are exactly the
{\it characters\/} of the group $S(\infty)$, discovered in \cite{T}. 

One of the main conclusions of \cite{VK} was that a solution to the
asymptotic problem mentioned above can be described in terms of the
limit behavior of the Frobenius coordinates of the diagrams $\nu$.
This explains why we find it so important to express $\chi^\nu_\sigma$
through the Frobenius coordinates of $\nu$ and not the conventional row
coordinates $\nu_1,\nu_2,\dots$. 

Theorem 1 is equivalent to the following claim. Let $\dim(\mu,\nu)$ be
the number of standard tableaux of skew 
shape $\nu/\mu$ provided that $\mu\subseteq\nu$, and 0 
otherwise.  

\proclaim{Theorem 2} For any Young diagram $\mu$ there exists an
inhomogeneous symmetric function 
$$
FS_\mu=s_\mu\,+\,\text{\rm lower degree terms}
$$
such that for any diagram $\nu=(p_1,\dots,p_d\mid q_1,\dots,q_d)$ 
with $|\nu|=n\ge m=|\mu|$,
$$ 
\frac{\dim(\mu,\nu)}{\dim\nu} 
=\frac{FS_\mu(p_1+\frac12,\dots,p_d+\frac12;q_1+\frac12,\dots,q_d+\frac12)} 
{n(n-1)\dots(n-m+1)}\,. 
$$
\endproclaim

Again, the function in question is defined uniquely. We call the
functions $FS_\mu$ the {\it Frobenius--Schur functions,\/} because
their definition exploits Frobenius coordinates.  

The expression $\frac{\dim(\mu,\nu)}{\dim\nu}$ has an important
meaning: it is a conditional probability in a 
probabilistic interpretation of the characters of the group
$S(\infty)$, which was suggested in \cite{VK}. We will briefly
explain this point, because it may help in understanding the role of
the Frobenius--Schur functions. 

Let $\Tab$ be the set of all infinite Young tableaux \cite{VK},
\cite{KV2}. An element $\tau\in\Tab$ is an infinite chain
$(\tau_0=\varnothing,\tau_1,\tau_2,\dots)$ 
of Young diagrams such that $|\tau_n|=n$ and $\tau_n$ is obtained
from $\tau_{n-1}$ by adding a box. In other terms, $\tau$ is an
infinite path in the Young graph. Given a probability measure $P$ on
$\Tab$, we can speak about random paths. We say that $P$ is a
{\it central measure\/} if the probability that the random path has a fixed
beginning $(\tau_0,\dots,\tau_n)$ depends only on the diagram
$\tau_n$ (and not on the whole finite chain $(\tau_0,\dots,\tau_n)$).
It turns out that there is a natural one--to--one
correspondence between the characters of $S(\infty)$ and
the indecomposable central measures on $\Tab$. Now we are in a position
to explain what $\frac{\dim(\mu,\nu)}{\dim\nu}$ is: this is the
probability that the random path passes through $\mu$ at ``time'' $m$
under the condition that it passes through $\nu$ at ``time'' $n$.
Note that this conditional probability does not depend on a central
measure; it is rather a characteristic of the Young graph. For more
detail, see \cite{VK} and \cite{KOr}

{}From quite a different point of view the Frobenius--Schur functions
are considered in \cite{BO, \S6}. 

\subhead The multiparameter Schur functions \endsubhead
Fix an arbitrary doubly infinite sequence $a=(a_i)_{i\in\Z}$ of
complex parameters. We define elements $h_{1;a}, h_{2;a},\dots$ of
$\La$ (which depend on the sequence $a$) by a generating series: 
$$ 
1+\sum_{k=1}^\infty \frac{h_{k;a}}{(u-a_1)\dots(u-a_k)} 
=1+\sum_{k=1}^\infty\frac{h_k}{u^k},
$$  
where $h_1,h_2,\ldots\in\La$ are the complete homogeneous
symmetric functions, $u$ is a formal variable, and the series above
are understood as elements of the algebra $\La[[\frac 1u]]$.
Equivalently, 
$$
h_{k;a}=\sum_{i=1}^k (-1)^{k-i}e_{k-i}(a_1,\dots,a_{k-1})\,h_i\,,
$$
where $e_1,e_2,\dots$ are the elementary symmetic functions. 
This implies, in particular, that $\{h_{k;a}\}_{k=1,2,\dots}$ is a 
system of algebraically independent generators of $\La$. 
 
We agree that $h_{0;a}=1$ and $h_{-1;a}=h_{-2;a}=\dots=0$. 
We also need the following notation: for $r\in\Z$, let $\tau^r a$ be the 
result of shifting $a$ by $r$ digits to the left:
$(\tau^r a)_i=a_{i+r}$.
 
The {\it multiparameter Schur function\/} indexed by an arbitrary
Young diagram $\mu$ is defined as follows:
$$ 
s_{\mu;a}=\det[h_{\mu_i-i+j; \tau^{1-j}a}],
$$ 
where the order of the determinant is any number greater or equal to 
$\ell(\mu)$, the number of rows in $\mu$. Clearly, 
$h_{k;a}=s_{(k);a}$. 

We have $s_{\mu;a}=s_\mu\,+\,\text{lower degree terms}$, which
implies that the functions $s_{\mu;a}$ form a basis in $\La$. When
$a\equiv0$, they reduce to $s_\mu$.

\proclaim{Theorem 3} The Frobenius--Schur functions $FS_\mu$ coincide
with the functions $s_{\mu;a}$ for the following special choice of the
parameters{\rm:} $a_i=i-\frac12$.
\endproclaim

Let us call $a_i=i-\frac12$ the {\it special parameter sequence.}

Although we are mainly interested in Frobenius--Schur functions, all
the results below are stated in a more general context of
multiparameter Schur functions. To pass to the Frobenius--Schur case
it suffices to substitute $a_i=i-\frac12$.  

\subhead Connections with earlier works \endsubhead 
Some topics in the 
representation theory of infinite--dimensional classical groups and
in the theory of Yangians naturally lead to the {\it
algebra of shifted symmetric functions.\/} This algebra, denoted by
$\La^*$, is studied in detail in \cite{OO}. This is a filtered algebra with
associated graded algebra canonically isomorphic to $\La$. In
$\La^*$, there is a natural basis formed by the {\it shifted Schur
functions.\/} These functions are closely related to the {\it factorial
Schur polynomials\/} discovered by Biedenharn and Louck \cite{BL1},
\cite{BL2} (for more references, see \cite{OO}). According to
\cite{KO}, both $\La^*$ and $\La$ can be realized as certain algebras of
functions on the set of partitions, which leads to an
algebra isomorphism $\La^*\to\La$. It turns out that {\it under this
isomorphism, the shifted Schur functions turn into the
Frobenius--Schur functions\/}: this claim is equivalent to Theorem
8.1 in \cite{OO} which expresses
$\frac{\dim(\mu,\nu)}{\dim\nu}$ in terms of the shifted 
Schur functions evaluated at the row coordinates of $\nu$ (not in
terms of the Frobenius coordinates!). A connection between 
shifted Schur functions and Schubert polynomials is explained in
\cite{L2}, see also Lascoux's remark in \cite{OO, \S15.7}. It is
possible that the Schubert polynomial techniques \cite{Ma3} can
provide an alternative approach to the results of the present paper
but we did not verify this.

Next, the functions $s_{\mu;a}$ can be obtained by an appropriate
specialization of Lascoux's {\it multi--Schur functions\/} introduced
in \cite{L1} (see also \cite{Ma3} and \cite{L2}). We are grateful to
Alain Lascoux for pointing out this fact.  For more detail see
\cite{ORV}.

Finally, the functions $s_{\mu;a}$ are closely connected with Molev's
{\it factorial supersymmetric Schur polynomials\/} \cite{Mo}. For
more detail see \cite{ORV}. At the early stage of the present work
this connection helped us very much. 

\subhead Duality \endsubhead
Consider the involutive automorphism $\om:\La\to\La$, see \cite{Ma1}.
Recall that 
$$
\om(\bp_k)=(-1)^{k-1}\bp_k,\quad 
\om(h_k)=e_k, \quad
\om(e_k)=h_k, \quad
\om(s_\mu)=s_{\mu'}\,.
$$
In terms of supersymmetric functions, $\om$ interchanges
the arguments $x,y\in\C^\infty_0$. 

Define the {\it dual parameter sequence\/} $\wha$ by
$\wha_i=-a_{-i+1}$. Thus, $a$ can be written as
$(\dots,-\wha_2,-\wha_1,a_1,a_2,\dots)$. 

\proclaim{Theorem 4} We have $\om(s_{\mu;a})=s_{\mu';\wha}$. 
\endproclaim

Set $e_{k;a}=s_{(1^k);a}$. By Theorem 4, $e_{k;a}=\om(h_{k;\wha})$
and we have a generating series
$$
1+\sum_{k=1}^\infty \frac{e_{k;a}}{(u-\wha_1)\dots(u-\wha_k)} 
=1+\sum_{k=1}^\infty\frac{e_k}{u^k}.
$$
Or, equivalenty, 
$$
e_{k;a}=\sum_{i=1}^k (-1)^{k-i}e_{k-i}(\wha_1,\dots,\wha_{k-1})\,e_i\,,
$$
Further, $s_{\mu;a}=\det[e_{\mu'_i-i+j; \tau^{j-1}a}]$.

These formulas can be immediately specialized for the case of the
Frobenius--Schur functions. Here 
$(a_1,a_2,\dots)=(\wha_1,\wha_2,\dots)=(\frac12,\frac32,\dots)$. Note
that the special parameter sequence is selfdual, so that
$\om(FS_\mu)=FS_{\mu'}$, which is also evident from Theorem 2.

\subhead Hook functions  \endsubhead
The definition of the multiparameter Schur functions has the form
$s_{\mu;a}=\det[h_{\mu_i-i+j, j-1}]$, where $\{h_{k,l}\}$ is a
certain family of elements of $\La$ (depending of $a$). By a general
result (see \cite{Ma2}, \cite{Ma1, Ex. I.3.21}), it follows

\proclaim{Theorem 5 (Giambelli formula)} Write $\mu$ in the Frobenius
notation, $\mu=(p_1,\dots,p_d\mid q_1,\dots,q_d)$. Then
$s_{\mu;a}=\det[s_{(p_i\mid q_j);a}]$.
\endproclaim 

The hook functions $s_{(p\mid q);a}$ can be determined from the
following generating series. Set 
$$
H(u)=1+\sum_{k=1}^\infty\frac{h_k}{u^k}\,, \qquad
E(v)=1+\sum_{k=1}^\infty\frac{e_k}{u^k}.
$$

\proclaim{Theorem 6} We have 
$$
1+(u+v)\sum_{p,q=0}^\infty
\frac{s_{(p\mid q);a}}
{(u-a_1)\dots(u-a_{p+1})(v-\wha_1)\dots(v-\wha_{q+1})}
=H(u)E(v).
$$
\endproclaim

When $a\equiv0$, this reduces to a known identity, see 
\cite{Ma1, Ex. I.3.9}. For the Frobenius--Schur functions we get
$$
1+(u+v)\sum_{p,q=0}^\infty
\frac{FS_{(p\mid q)}}
{(u-\frac12)\dots(u-\frac{2p+1}2)(v-\frac12)\dots(v-\frac{2q+1}2)}
=H(u)E(v).
$$

\subhead Transition coefficients \endsubhead
Let $b=(b_i)_{i\in\Z}$ be another parameter sequence. By Theorem 6,
$$
\gather
1+(u+v)\sum_{p,q=0}^\infty
\frac{s_{(p\mid q);a}}
{(u-a_1)\dots(u-a_{p+1})(v-\wha_1)\dots(v-\wha_{q+1})}\\
=1+(u+v)\sum_{p,q=0}^\infty
\frac{s_{(p\mid q);b}}
{(u-b_1)\dots(u-b_{p+1})(v-\whb_1)\dots(v-\whb_{q+1})}\,.
\endgather
$$
It follows that
$$ 
s_{(p\,|\,q);a}=\sum_{p'=0}^p\sum_{q'=0}^q  
c_{pp'}(a,b)\,c_{qq'}(\wha,\whb) \,s_{(p'\,|\,q');b}\,,
$$
where
$$ 
c_{pp'}(a,b)=h_{p-p'}(b_1,\dots,b_{p'+1};-a_1,\dots,-a_p), 
\qquad p\ge p'\ge0.
$$
Note that $c_{pp}(a,b)=1$.

An important particular case is $b\equiv0$. Then
$$ 
c_{pp'}(a,0)=h_{p-p'}(0,\dots,0;-a_1,\dots,-a_p) 
=(-1)^{p-p'}e_{p-p'}(a_1,\dots,a_p).
$$

Together with the Giambelli formula (Theorem 5) these formulas
provide a simple determinant expression for the transition coefficients
between the bases $\{s_{\mu;a}\}$ and $\{s_{\mu;b}\}$. 

\proclaim{Theorem 7} Let $a=(a_i)_{i\in\Z}$ and 
$b=(b_i)_{i\in\Z}$ 
be two sequences of parameters. In the expansion 
$$ 
s_{\mu;a}=\sum_\nu c_{\mu\nu}(a,b)\,s_{\nu;b} 
$$ 
the coefficients $c_{\mu\nu}(a,b)$ vanish unless $\nu\subseteq\mu$ 
and both diagrams have the same number of diagonal squares{\rm:}
$d(\nu)=d(\mu)$.  
 
Assume $d(\nu)=d(\mu)=d$ and write both diagrams in the Frobenius 
notation,  
$$ 
\mu=(p_1,\dots,p_d\,|\,q_1,\dots,q_d),\quad 
\nu=(p'_1,\dots,p'_d\,|\,q'_1,\dots,q'_d). 
$$ 
Then   
$$ 
c_{\mu\nu}(a,b)=\det[c_{p_i,p'_j}(a,b)]\, 
\det[c_{q_i,q'_j}(\wha,\whb)], 
$$ 
where the determinants are of order $d$ and the coefficients of the
form $c_{pp'}(a,b)$ are defined above. 
\endproclaim 

In particular, taking $b=0$ we get a nice expression of the
multiparameter Schur functions
$s_{\mu;a}$ through the conventional Schur functions. From this
expression we get

\proclaim{Corollary 8} Let $p=\mu_1-1$ and $q=\mu'_1-1$. The function
$s_{\mu;a}$ depends only of $p+q$ parameters $(a_1,\dots,a_p)$ and
$(\wha_1,\dots,\wha_q)=(-a_0,\dots,-a_{q-1})$. 
\endproclaim 

Let us specially state the particular case of Theorem 7 corresponding
to $a_i=i-\frac12$ and $b\equiv0$: 

\proclaim{Theorem 9} We have
$$
FS_{(p_1,\dots,p_d\mid q_1,\dots,q_d)}=
\sum \Sb p'_1\le p_1, \dots, p'_d\le p_d\\
q'_1\le q_1,\dots, q'_d\le q_d \endSb
\det[c_{p_i,p'_j}]\,\det[c_{q_i,q'_j}]\,
s_{(p'_1,\dots,p'_d \mid q'_1,\dots,q'_d)}\,,
$$
where
$$
c_{pp'}=(-1)^{p-p'}e_{p-p'}(\underbrace{\tfrac12,\dots,\tfrac{2p-1}2}_p).
$$
\endproclaim

\subhead Vanishing, interpolation, characterization \endsubhead
For a nonempty Young diagram written in the Frobenius notation
$$ 
\la=(p_1(\la),\dots,p_d(\la)\mid q_1(\la),\dots,q_d(\la)),
\qquad d=d(\la),  
$$ 
we set 
$$  
(x(\la);y(\la))=
(a_{p_1(\la)+1},\dots,a_{p_d(\la)+1};
\wha_{q_1(\la)+1},\dots,\wha_{q_d(\la)+1})\in
\C^d\times\C^d\subset\C^\infty_0\times\C^\infty_0.  
$$ 
Note that for the special parameter sequence $a_i=i-\half$, we have
$$
(x(\la);y(\la))=(p_1(\la)+\half,\dots,p_d(\la)+\half;
q_1(\la)+\half,\dots,q_d(\la)+\half).
$$
We also agree that $(x(\varnothing);y(\varnothing))=(0;0)$.
 
\proclaim{Theorem 10} If 
$\mu\not\subseteq\la$ then $s_{\mu;a}(x(\la);y(\la))=0$. 
\endproclaim 

\proclaim{Theorem 11} We have 
$$
s_{\mu;a}(x(\mu);y(\mu))=
\prod_{(i,j)\in\mu}(a_{\mu_i-i+1}-a_{j-\mu'_j})
=\prod_{(i,j)\in\mu}(a_{\mu_i-i+1}+\wha_{\mu'_j-j+1}).
$$
If the numbers $a_i$ are pairwise distinct then
$s_{\mu;a}(x(\mu);y(\mu))\ne0$ for all $\mu$.
\endproclaim

For the Frobenius--Schur functions this reduces to 
$$
FS_\mu(p_1(\mu)+\half,\dots,p_d(\mu)+\half;
q_1(\mu)+\half,\dots,q_d(\mu)+\half)=
\prod_{(i,j)\in\mu}(\mu_i-j+\mu'_j-i+1).
$$
By the hook formula, the last expression equals $\dim\mu/|\mu|!$,
which agrees with Theorem 2.

As a corollary of Theorem 10 and Theorem 11 we get

\proclaim{Theorem 12} Assume that the parameters $a_i$ are pairwise
distinct. Given an arbitrary element $f\in\La$ of degree $\le n$,
consider its expansion in the multiparameter Schur functions of
degree $\le n${\rm:}
$f=\sum\limits_{\la:\, |\la|\le n}c(\la)s_{\la;a}$.  
Then the coefficients $c(\la)$ can be found by recurrence on $|\la|$,
as follows{\rm:} 
$$
c(\varnothing)=f(0;0), \qquad
c(\la)=\frac{f(x(\la);y(\la))
-\sum\limits_{\mu:\,\mu\subset\la, 
\mu\ne\la}c(\mu)s_{\mu;a}(x(\la);y(\la))}
{s_{\la;a}(x(\la);y(\la))}\,.
$$
\endproclaim

In other words, the multiparameter Schur functions naturally arise in
the following Newton--type {\it interpolation process\/}: we consider the
space of supersymmetric functions on $\C^\infty_0\times\C^\infty_0$ with its
natural filtration by degree and take the points of interpolation of
the form $(0;0)$ and
$(a_{i_1},\dots,a_{i_d};\wha_{j_1},\dots,\wha_{j_d})$. Here 
$d=1,2,\dots$; $i_1>\dots>i_d\ge1$ and $j_1>\dots>j_d\ge1$ are 
$d$-tuples of natural numbers; $(a_1,a_2,\dots)$ and $(\wha_1,\wha_2,\dots)$
are two fixed sequences of complex parameters such that $a_i\ne a_j$ and
$\wha_i\ne\wha_j$ for $i\ne j$, and $a_i\ne-\wha_j$ for all $i,j$ (the
condition $a_i\ne-\wha_j$ is natural because of the cancellation
property).  

Note that the Frobenius--Schur functions correspond to the
interpolation process with the points of interpolation of the form
$(0;0)$ and $(x_1>\dots>x_d;y_1>\dots>y_d)$, where $x_i, y_j$ are
taken from the set $\{\frac12,\frac32,\frac52,\dots\}$.

\proclaim{Theorem 13} Assume that the parameters $a_i$ are pairwise
distinct. Then the function $F=s_{\mu;a}\in\La$ can be characterized by the
following properties {\rm(i)} and {\rm(ii)} as well as by the
properties {\rm(i${}'$)} and {\rm(ii${}'$):}

{\rm(i)} $F$ vanishes at $(x(\la);y(\la))$ for all diagrams
$\la\ne\mu$ such that $|\la|\le|\mu|$;

{\rm(ii)} $F$ is of degree $\le|\mu|$ and its value at 
$(x(\mu);y(\mu))$ is the same as in Theorem 11.

{\rm(i${}'$)} $F$ vanishes at $(x(\la);y(\la))$ for all diagrams
$\la$ with $|\la|<|\mu|$;

{\rm(ii${}'$)} the top degree homogeneous component of $F$ coincides with
the Schur function $s_\mu$.
\endproclaim

In the particular case of the Frobenius--Schur functions, the results
of this section are equivalent to similar claims for the shifted
Schur functions contained in \cite{Ok1}, \cite{OO, \S3}.
Interpolation of arbitrary polynomials in terms of Schubert
polynomials is discussed in \cite{L2}. A very general scheme of
Newton interpolation for symmetric polynomials is developed in
\cite{Ok2}. In the latter paper one can also find the references to earlier
works by Knop, Sahi, and Okounkov.

\subhead The combinatorial formula \endsubhead
Let $\Z'=\{\dots, -\tfrac32, -\tfrac12, \tfrac12,\tfrac32,\dots\}$ 
stand for the set of proper half--integers. For any $\ep\in\Z'$ set
$a'_\ep=a_{\ep+1/2}$. Note that $(\wha)'_\ep=a'_{-\ep}$. For the
special parameter sequence, $a'_\ep=\ep$. 

Let $\nu$ be a skew hook (in another terminology, a ribbon). That is,
$\nu$ is a connected skew Young diagram containing no $2\times 2$
block of squares. We attach to $\nu$ a polynomial $f_{\nu;a}(u,v)$ in
two variables $u,v$, of degree $|\nu|$, as follows. 

Consider the interior sides $s$ 
of $\nu$: each interior side is adjacent to two squares of 
$\nu$; the total number of the interior sides is equal to $|\nu|-1$. 
To each interior side $s$ we attach the coordinates $(\ep,\de)$ of 
its midpoint, \footnote{As in \cite{Ma1}, we assume that the first
coordinate axis is directed downwards and the second coordinate axis
is directed to the right.} and we write $s=(\ep,\de)$. Note that one of the  
coordinates is always in $\Z'$  while another coordinate is 
integral. Specifically, if $s$ is a vertical side then $\ep\in\Z'$, 
$\de\in\Z$, and the ends of $s$ are the points $(\ep-1/2,\de)$ and 
$(\ep+1/2,\de)$; if $s$ is a horizontal side then $\ep\in\Z$, 
$\de\in\Z'$, and the ends of $s$ are the points $(\ep,\de-1/2)$, 
$(\ep,\de+1/2)$.  

For both vertical and horizontal sides, $\de-\ep\in\Z'$, so that
$a'_{\de-\ep}=a_{\de-\ep+1/2}$ is well defined. We set  
$$ 
f_{\nu;a}(u,v)=(u+v) 
\prod\Sb \text{vertical interior}\\ \text{sides $s=(\ep,\de)$ of $\nu$} 
\endSb (u-a'_{\de-\ep}) 
\prod\Sb \text{horizontal interior}\\ \text{sides $s=(\ep,\de)$ of $\nu$} 
\endSb (v+a'_{\de-\ep}). 
$$ 
 
For instance, if $\nu=(4,2,2)/(1,1)$ then there are 6 squares and 5 
interior sides with midpoints
$$ 
(\tfrac52,1),\quad(2,\tfrac32),\quad(1,\tfrac32),\quad 
(\tfrac12,2),\quad(\tfrac12,3), 
$$ 
and we have  
$$ 
\gather 
f_{\nu;a}(u,v)=(u+v)(u-a'_{-3/2})(u-a'_{3/2})(u-a'_{5/2}) 
(v+a'_{-1/2})(v+a'_{1/2})\\ 
=(u+v)(u-a_{-1})(u-a_2)(u-a_3)(v+a_0)(v+a_1) 
\endgather 
$$

Note that $f_{\nu;a}(v,u)=f_{\nu';\wha}(u,v)$. 

More generally, if $\nu$ is an arbitrary skew diagram with no
$2\times2$ block of squares, then $\nu$ is a disjoint union of skew
hooks and we define $f_{\nu;a}(u,v)$ as the product of the polynomials 
attached to the connected components of $\nu$. 

We also agree that $f_\varnothing(u,v)\equiv1$. Note that if
$\nu\ne\varnothing$ then $f_{\nu;a}(0;0)=0$, because of the factor $u+v$.

Given a Young diagram $\mu$, a {\it diagonal--strict\/} tableau of
shape $\mu$ is a function 
$\T(i,j)$ from the squares $(i,j)$ of $\mu$ to  
$\{1,2,\dots\}$ such that the numbers $\T(i,j)$ weakly increase 
both along the rows $i=\const$ (from left to right) and down the
columns $j=\const$, and 
strictly increase along the diagonals $j-i=\const$. Each subset of
form $\T^{-1}(k)$ is a skew diagram with no $2\times2$ block of
squares, so that the polynomial $f_{\T^{-1}(k)}(u,v)$ makes sense. 

\proclaim{Theorem 14} The following combinatorial formula holds
$$
s_{\mu;a}(x;y)=\sum_\T 
\prod_{k\ge1}f_{\T^{-1}(k);a}(x_k,y_k),
$$ 
summed over all diagonal--strict tableaux of shape $\mu$. 
\endproclaim

If $n$ is so large that
$x_{n+1}=x_{n+2}=\dots=y_{n+1}=y_{n+2}=\dots=0$ then only tableaux
$\T$ with values in $\{1,\dots,n\}$ really contribute. Further, 
each product in the above formula is actually finite, because
$\T^{-1}(k)=\varnothing$ for sufficiently large $k$ (and then
$f_{\T^{-1}(k)}\equiv1$). 
Thus, the right--hand side of the combinatorial formula makes sense.

The formula is compatible with the involution
$\om$, because the definition of a diagonal--strict tableau
is symmetric with respect to transposition and
$f_{\nu;a}(v,u)=f_{\nu';\wha}(u,v)$. 

The formula can be restated in terms of tableaux of a different kind,
see \cite{ORV}.
The proof of Theorem 14 is due to Vladimir Ivanov. 

Note that for the first time, an ``inhomogeneous'' combinatorial
formula probably appeared in \cite{BL1}, \cite{BL2}, see also
\cite{CL}. Other examples can be found in \cite{GG},
\cite{Ma2}, \cite{Mo}, \cite{Ok1}, \cite{OO}, \cite{Ok2}. See also
further references in \cite{Ok2} to works by Knop, Okounkov, and Sahi
about combinatorial formulas for interpolation Jack and Macdonald
polynomials.  

\subhead The Sergeev--Pragacz formula \endsubhead
In this section we assume that $(x;y)$ ranges over
$\C^n\times\C^n\subset\C^\infty_0\times\C^\infty_0$. Thus, we deal
with supersymmetric polynomials in $n+n$ variables. Then
$s_{\mu;a}(x;y)$ vanishes unless $d(\mu)\le n$.

In the space of all polynomials in $n+n$ variables there is a natural
action of $\frak S_n\times\frak S_n$, the product of two
copies of the symmetric group $\frak S_n$. For 
$w=(w_1,w_2)\in\frak S_n\times\frak S_n$ we write
$\ep(w)=\sgn(w_1)\sgn(w_2)$. Set 
$$
\gather 
V(x_1,\dots,x_m)=\prod_{1\le i<k\le m}(x_i-x_k), \qquad 
V(y_1,\dots,y_n)=\prod_{1\le j<l\le n}(y_j-y_l), \\
(x\,|\,a)^m=\cases (x-a_1)\dots(x-a_m), & m\ge1, \\ 
1, & m=0, \endcases
\endgather 
$$
and $(k)_+=\max(k,0)$.

\proclaim{Theorem 15} Fix $n=1,2,\dots$ and let $\mu$ be an arbitrary 
Young diagram such that $d=d(\mu)\le n$. We have 
$$ 
s_{\mu;a}(x_1,\dots,x_n;y_1,\dots,y_n)= 
\frac{\sum\limits_{w\in\frak S_n\times\frak S_n}\ep(w)\,  
w[g_{\mu;a}(x_1,\dots,x_n;y_1,\dots,y_n)]} 
{V(x_1,\dots,x_n)V(y_1,\dots,y_n)}\,,  
$$ 
where 
$$ 
\gather 
g_{\mu;a}(x_1,\dots,x_n;y_1,\dots,y_n)= 
\prod_{i=1}^d(x_i\,|\,a)^{\mu_i-i}x_i^{(n-\mu_i)_+} 
(y_i\,|\,\wha)^{\mu'_i-i}y_i^{(n-\mu'_i)_+}\\ 
\times\prod_{i=d+1}^n x_i^{n-i} y_i^{n-i}\cdot 
\prod\Sb i,j\le n\\ (i,j)\in\mu\endSb (x_i+y_j). 
\endgather 
$$
\endproclaim  

When $a\equiv0$, this turns into the conventional Sergeev--Pragacz
formula, see \cite{PT}, \cite{Ma3}, \cite{Ma1, Ex. I.3.23}. Our proof
is an adaptation of the argument of \cite{PT}. 

Taking $n=d$ we get an analog of the factorization property
\cite{BR}: 

\proclaim{Corollary 16} Let $\mu=(p_1,\dots,p_d\,|\,q_1,\dots,q_d)$ 
be a Young diagram of depth $d$. Then 
$$ 
s_{\mu;a}(x_1,\dots,x_d;y_1,\dots,y_d)= 
\frac{\det[(x_i\,|\,a)^{p_j}]_{i,j=1}^d}{V(x_1,\dots,x_d)}\cdot 
\frac{\det[(y_i\,|\,\wha)^{q_j}]_{i,j=1}^d}{V(y_1,\dots,y_d)}\cdot 
\prod_{i,j=1}^d (x_i+y_j)\,. 
$$ 
\endproclaim

\subhead Schur's $Q$-functions \endsubhead
As is well known, there exists a deep analogy between the
conventional Schur functions (also called Schur's $S$-functions) and
Schur's $Q$-functions. The latter are related to projective
representations of symmetric groups just in the same way as the former
are related to ordinary representations. Schur's $Q$-functions span
a subalgebra in $\La$, whose elements can be characterized by a
supersymmetry property of another kind, see, e.g., \cite{Ma1, \S
III.8}, \cite{P}. The results of the present paper have counteparts
for Schur's $Q$-functions: these are due to Vladimir Ivanov (\cite{I1},
\cite{I2}, and a paper in preparation).

\bigskip

{\smc G.~Olshanski}: Dobrushin Mathematics Laboratory, Institute for
Information Transmission Problems, Bolshoy Karetny 19, 
Moscow 101447, GSP-4, Russia.  

E-mail address: {\tt olsh\@iitp.ru, olsh\@glasnet.ru}

{\smc A.~Regev}: Department of Theoretical Mathematics,  
Weizmann Institute of Science, Rehovot 76100, Israel.

E-mail address: {\tt regev\@wisdom.weizmann.ac.il} 

{\smc A.~Vershik}: Steklov Mathematical Institute (POMI), Fontanka 27, 
St.~Petersburg 191011, Russia.

E-mail address: {\tt vershik\@pdmi.ras.ru}

\enddocument

\bye